\pdfoutput=1
\documentclass[11pt]{article}
\usepackage[margin=1in]{geometry}
\usepackage[T1]{fontenc}
\usepackage{lmodern,microtype}
\usepackage{amsmath,amssymb,amsthm}
\usepackage{booktabs}
\usepackage{xcolor}
\usepackage{tikz}
\usetikzlibrary{arrows.meta,shapes.misc,decorations.pathmorphing,decorations.pathreplacing,patterns,calc}
\usepackage[colorlinks=true,linkcolor=blue!50!black,citecolor=blue!50!black,urlcolor=blue!50!black,bookmarksnumbered=true]{hyperref}
\usepackage{bookmark}
\hypersetup{pdftitle={No infinite cluster at criticality for long-range percolation with fast-decaying interactions},pdfauthor={Yan Ru Pei},pdfkeywords={long-range percolation, critical percolation, continuity of the phase transition, gluing inequality, renormalization}}
\newtheorem{theorem}{Theorem}[section]

\newtheorem{lemma}[theorem]{Lemma}
\newtheorem{corollary}[theorem]{Corollary}
\theoremstyle{definition}
\newtheorem{definition}[theorem]{Definition}
\theoremstyle{remark}
\newtheorem{remark}[theorem]{Remark}
\numberwithin{equation}{section}
\newcommand{\Z}{\mathbb Z}
\newcommand{\R}{\mathbb R}

\newcommand{\Pp}{\mathbb P}
\newcommand{\E}{\mathbb E}
\newcommand{\1}{\mathbf 1}
\newcommand{\cC}{\mathcal C}
\newcommand{\cI}{\mathcal I}
\newcommand{\cJ}{\mathcal J}
\newcommand{\cL}{\mathcal L}
\newcommand{\lr}{\leftrightarrow}
\newcommand{\nlr}{\mathrel{\not\leftrightarrow}}
\newcommand{\Lam}{\Lambda}
\newcommand{\norm}[1]{\|#1\|}
\makeatletter\g@addto@macro\bfseries{\boldmath}\makeatother

\title{No infinite cluster at criticality for long-range percolation\\ with fast-decaying interactions}
\author{Yan Ru Pei\\[0.3em]\small\href{mailto:yanrpei@gmail.com}{\texttt{yanrpei@gmail.com}}}
\date{October 2026}

\begin{document}
\maketitle

\begin{abstract}
We consider independent long-range bond percolation on $\mathbb{Z}^d$ in which $x$ and $y$ are
joined with probability $1-e^{-\beta J(x-y)}$. Let $d\ge2$ and let $J$ be invariant under coordinate
permutations and sign changes, positive on nearest-neighbour pairs, and satisfy
$\liminf_{n\to\infty} n^d\sum_{\|x\|>n}J(x)=0$. We prove that there is no infinite cluster at the
critical parameter. Together with Berger's theorem for heavy tails, this gives continuity of the
phase transition for $J(x)\asymp\|x\|^{-d-\alpha}$ for every $\alpha>0$ except $\alpha=d$, in every
dimension $d\ge2$, and for every exponentially decaying kernel satisfying these hypotheses. Previously the
fast-decaying case was known for nearest-neighbour and uniform finite-range kernels, in $d=2$ for
finite-range and exponential kernels, and in high dimension for sufficiently spread-out kernels.
The proof adapts a recent renormalization argument for nearest-neighbour percolation on
$\mathbb{Z}^3$ that is driven by a gluing inequality conjectured by Kozma and Nitzan. Long edges enter in two
places: they can jump over the surfaces that a nearest-neighbour path must cross, and they let a
cluster outside a box enter it from far away. We control the first with thick outer faces and the
tail condition, and the second by measuring entrances through the total edge weight between a
cluster and a box rather than through the number of boundary vertices.
\end{abstract}

\section{Introduction}\label{sec:intro}

\subsection{Long-range percolation and the main result}

Let $d\ge1$ and let $J:\Z^d\setminus\{0\}\to[0,\infty)$ be a kernel with $J(x)=J(-x)$ and
$\sum_xJ(x)<\infty$. In
\emph{long-range percolation} at parameter $\beta>0$, each unordered pair $\{x,y\}$ of distinct
vertices of $\Z^d$ is open independently with probability
\[
p_{xy}(\beta)=1-e^{-\beta J(x-y)},
\]
and closed otherwise; pairs with $J(x-y)=0$ are never open. Write $\Pp_\beta$ for the law,
$\theta(\beta)=\Pp_\beta(0\lr\infty)$ for the probability that the origin lies in an infinite open
cluster, and $\beta_c=\inf\{\beta:\theta(\beta)>0\}$. The function $\theta$ is nondecreasing, and
the basic question is whether it is continuous at $\beta_c$, that is, whether $\theta(\beta_c)=0$.

For nearest-neighbour percolation, $\theta(p_c)=0$ was proved in $d=2$ by Harris and
Kesten~\cite{Harris1960,Kesten1980} and in high dimensions by the lace
expansion~\cite{HaraSlade1990,FitznerHofstad2017}. Kozma and Nitzan~\cite{KozmaNitzan2024} reduced
the remaining dimensions to a correlation inequality for independent percolation on finite graphs,
which has since been proved in several forms~\cite{Leder2026,OpenAIZ3,BlancRenaudie2026,Gouere2026}. We use it in the
joint form of Theorem~\ref{thm:gluing}, which holds for arbitrary edge probabilities; all graphs to
which we apply it are finite.

For long-range models the answer depends on the tail of $J$. Writing
\[
T(n)=\sum_{\norm{x}>n}J(x),\qquad \norm{x}=\max_i|x_i|,
\]
we use the following hypotheses.
\begin{itemize}
\item[(Sym)] $J$ is invariant under coordinate permutations and sign changes of coordinates.
\item[(Pos)] $J(e_1)>0$, where $e_1$ is the first unit vector.
\item[(Tail)] $\liminf_{n\to\infty}n^dT(n)=0$.
\end{itemize}

\begin{theorem}\label{thm:main}
Let $d\ge2$ and let $J$ satisfy \textup{(Sym)}, \textup{(Pos)} and \textup{(Tail)}. Then
$0<\beta_c<\infty$ and $\Pp_{\beta_c}(0\lr\infty)=0$.
\end{theorem}

Combined with Berger's theorem for heavy-tailed kernels~\cite[Theorem~2.6]{Berger2002}, which we
recall in Section~\ref{sec:corollaries}, this settles the power-law case away from the borderline
exponent.

\begin{corollary}\label{cor:power}
Let $d\ge2$, $\alpha>0$, $\alpha\ne d$, and let $J$ satisfy \textup{(Sym)} and
$c\norm{x}^{-d-\alpha}\le J(x)\le C\norm{x}^{-d-\alpha}$ for all $x\ne0$, with $0<c\le C<\infty$.
Then $\theta(\beta_c)=0$.
\end{corollary}

\begin{corollary}\label{cor:exp}
Let $d\ge2$ and let $J$ satisfy \textup{(Sym)}, \textup{(Pos)} and $J(x)\le Ae^{-a\norm{x}}$ for some
$A,a>0$. Then $\theta(\beta_c)=0$. The same holds for every finite-range $J$ satisfying
\textup{(Sym)} and \textup{(Pos)}.
\end{corollary}

Corollary~\ref{cor:power} covers both regimes of $\alpha$ around the exceptional value $d$. For
$\alpha<d$ it is Berger's theorem; for $\alpha>d$ it follows from Theorem~\ref{thm:main}, since
then $T(n)\le C'n^{-\alpha}$. The value $\alpha=d$ is excluded: there $n^dT(n)$ is bounded away
from zero, so edges of length comparable to $n$ leave the cube $[-n,n]^d$ with probability of
order one at every scale. This is the mechanism behind the discontinuous transition of the
one-dimensional model with $J(x)\asymp|x|^{-2}$~\cite{AizenmanNewman1986,DuminilCopinGarbanTassion2020},
and Theorem~\ref{thm:main} says nothing about it.

\subsection{Earlier results}

Berger~\cite[Theorem~2.6]{Berger2002} proved that $\theta(\beta_c)=0$ whenever
$\liminf_{\norm{k}\to\infty}\Pp(0\lr k\text{ directly})\,\norm{k}_1^{s}>0$ for some $s<2d$, in every
$d\ge1$; Hutchcroft~\cite{Hutchcroft2021} made this quantitative. In the opposite regime the
kernel decays fast enough that the model should behave like nearest-neighbour percolation, and
Hutchcroft writes that ``it is widely believed that the phase transition should be continuous for
all $\alpha>0$ and $d\ge2$''~\cite[p.~2]{Hutchcroft2021}. The known cases of this regime are the
following.
\begin{itemize}
\item In $d=2$, Duminil-Copin, Sidoravicius and Tassion~\cite{DuminilCopinSidoraviciusTassion2016}
proved continuity for finite-range kernels invariant under the symmetries of $\Z^2$, and
Angel, de la Riva, Panagiotis and Shi~\cite{AngelEtAl2026} extended this to kernels with
exponential decay; they state the case $\alpha>2$ of power laws in $d=2$ as a hope
\cite[\S1.4]{AngelEtAl2026}. Both arguments are planar in spirit: they glue clusters that cross
as subsets of $\R^2$.
\item In high dimensions, the lace expansion gives mean-field behaviour, and in particular
$\theta(\beta_c)=0$, for sufficiently spread-out kernels when $d>3(\alpha\wedge2)$
\cite{HeydenreichHofstadSakai2008}.
\item For kernels that are uniform on a finite range, Theorem~\ref{thm:main}'s conclusion follows
from the recent preprint on quasi-transitive graphs~\cite{OpenAIQT}.
\end{itemize}
Theorem~\ref{thm:main} treats every $d\ge2$ and requires neither planarity nor high dimension
nor a spread-out kernel. It contains the theorem of~\cite{AngelEtAl2026}, whose kernels satisfy
(Sym) (in $d=2$, invariance under $\pi/2$-rotations and reflections), (Pos) and exponential decay.

\subsection{Idea of the proof}\label{sec:ideas}

\begin{figure}[tp]
\centering
\newcommand*{\rmloc}[1]{\\[1pt]{\tiny\color{black!75}#1}}%
\begin{NoHyper}%
\begin{tikzpicture}[
  x=1cm, y=1cm, font=\scriptsize,
  execute at begin node={\thickmuskip=5mu\medmuskip=4mu\relax},
  >={Stealth[length=1.5mm,width=1.1mm]},
  thmlook/.style={fill=black!80, draw=black, line width=0.8pt},
  ourslook/.style={rounded corners=1.3mm, fill=black!16, draw=black, line width=0.6pt},
  knownlook/.style={fill=white, draw=black, line width=0.4pt},
  redolook/.style={fill=white, draw=black, line width=0.45pt,
    path picture={\draw[black, line width=0.45pt]
      ([shift={(1.6pt,1.6pt)}]path picture bounding box.south west) rectangle
      ([shift={(-1.6pt,-1.6pt)}]path picture bounding box.north east);}},
  toollook/.style={rounded rectangle, fill=black!9, draw=black, line width=0.6pt},
  box/.style={rectangle, align=center, text width=2.75cm, inner xsep=1.2pt, inner ysep=1.5pt,
    minimum height=0.9cm},
  wide/.style={text width=4.37cm},
  thm/.style={box, wide, thmlook, text=white, font=\scriptsize\bfseries},
  ours/.style={box, ourslook},
  known/.style={box, knownlook},
  redo/.style={box, redolook, inner xsep=2.6pt, inner ysep=2.8pt},
  tool/.style={toollook, align=center, inner sep=2.5pt, minimum height=0.5cm},
  sw/.style={rectangle, minimum width=0.5cm, minimum height=0.3cm, inner sep=0pt},
  lg/.style={anchor=west, inner sep=2pt, text height=1.6ex, text depth=0.3ex},
  use/.style={->, line width=0.5pt, draw=black!80},
  rel/.style={use, dashed},
  band/.style={->, line width=0.9pt, draw=black!45},
  col/.style={rounded corners=2.5mm, fill=black!3, draw=black!30, line width=0.4pt},
]
\draw[col] (0.00,0.85) rectangle (9.90,-7.05);
\draw[col] (10.25,0.85) rectangle (16.45,-3.05);
% theorems
\node[thm] (T1) at (4.95,0) {Theorem~\ref{thm:main}\\ no infinite cluster at $\beta_c$\\ under (Sym), (Pos), (Tail)};
\node[thm] (T2) at (13.35,0) {Corollaries~\ref{cor:power}, \ref{cor:exp}\\ power laws, $\alpha\neq d$;\\ exponential tails};
% column 1
\node[redo, wide] (E6) at (4.95,-2.0) {Adaptive exploration\\ and contour count\rmloc{\cite[\S6]{OpenAIZ3}; here Sec.~\ref{sec:exploration}}};
\node[redo]  (R5) at (2.10,-4.0) {Relay between\\ coarse cubes\rmloc{\cite[Lem.~5.1]{OpenAIZ3}; here Lem.~\ref{lem:relay}}};
\node[redo]  (F3) at (7.80,-4.0) {Outer thick faces\\ at $\beta$\rmloc{\cite[Lem.~3.3]{OpenAIZ3}; here Lem.~\ref{lem:faces}}};
\node[redo]  (X4) at (2.10,-6.0) {Extension through\\ a fresh region\rmloc{\cite[Lem.~4.1]{OpenAIZ3}; here Lem.~\ref{lem:extension}}};
\node[known] (G2) at (7.80,-6.0) {Gluing inequality\rmloc{\cite[Prop.~2.1]{OpenAIZ3}, \cite[Thm~1]{Gouere2026}}};
\draw[use] (E6) -- (E6 |- T1.south);
\draw[use] (R5.north) -- (E6.south -| 3.45,0);
\draw[use] (F3.north) -- (E6.south -| 6.45,0);
\draw[use] (X4) -- (R5);
\draw[use] (G2) -- (X4);
\draw[use] (F3.south west) -- (X4.north east);
% column 2
\node[known, text width=2.5cm] (B1) at (11.75,-2.0) {Berger's theorem,\\ $s<2d$\rmloc{\cite[Thm~2.6]{Berger2002}}};
\node[known, text width=2.5cm] (A1) at (14.95,-2.0) {$d=2$, exponential\\ decay\rmloc{\cite[Thm~1.1]{AngelEtAl2026}}};
\draw[use] (B1) -- (B1 |- T2.south);
\draw[rel] (A1) -- (A1 |- T2.south);
\draw[use] (T1) -- (T2);
% band
\node[tool, text width=9.3cm] (S1) at (4.95,-7.75) {Scales $m\le w\ll R\ll r$ chosen from the tail condition (Lemma~\ref{lem:scales})};
\draw[band] (2.10,-7.75 |- S1.north) -- (2.10,-7.05);
\draw[band] (7.80,-7.75 |- S1.north) -- (7.80,-7.05);
% legend
\node[sw, thmlook]  (g1) at (0.45,-8.55) {}; \node[lg] at (g1.east) {main theorem};
\node[sw, toollook, minimum width=0.8cm] (g3) at (3.35,-8.55) {}; \node[lg] at (g3.east) {shared technique};
\node[sw, knownlook] (g4) at (0.45,-9.05) {}; \node[lg] at (g4.east) {existing result};
\node[sw, redolook]  (g5) at (3.35,-9.05) {}; \node[lg] at (g5.east) {existing result, adapted to long edges and to $\mathbb Z^d$};
\draw[use]  (0.20,-9.55) -- (0.70,-9.55) node[lg] {used in proof};
\draw[rel]  (3.10,-9.55) -- (3.60,-9.55) node[lg] {other relation};
\draw[band] (6.05,-9.55) -- (6.85,-9.55) node[lg] {technique used in the column};
\end{tikzpicture}%
\end{NoHyper}
\caption{How Theorem~\ref{thm:main} and its corollaries are proved. A solid arrow means that its source
is used in the proof of its target; only the main dependencies are drawn. The dashed arrow means that
the theorem of~\cite{AngelEtAl2026} is contained in the case $d=2$ of Corollary~\ref{cor:exp}. Box styles
give provenance; grey labels give the source of each adapted step and where it is proved here. The
modifications are listed in Section~\ref{sec:ideas}.}
\label{fig:roadmap}
\end{figure}
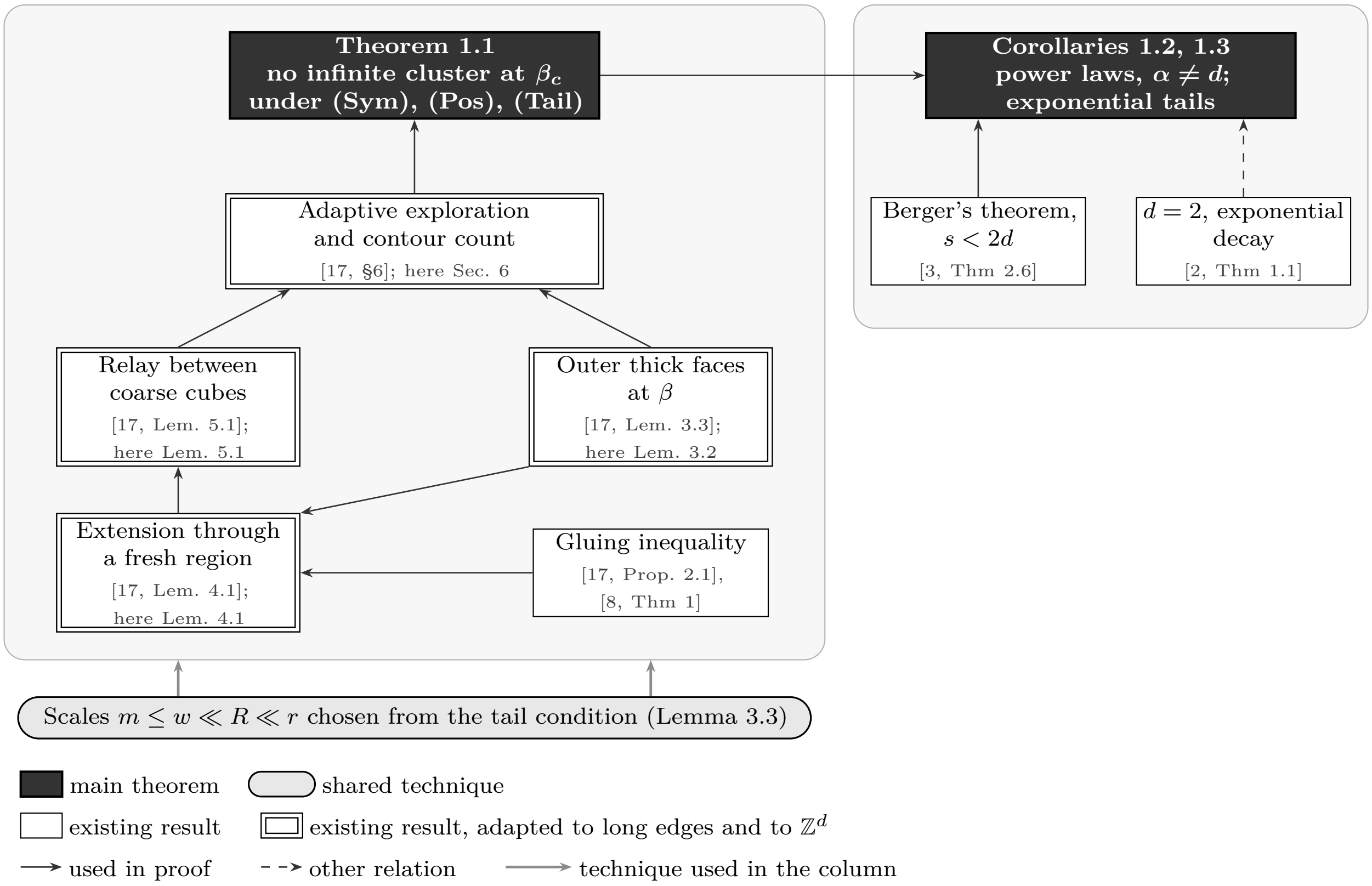

We follow the self-contained proof of $\theta(p_c)=0$ for nearest-neighbour percolation on $\Z^3$
in~\cite{OpenAIZ3}, which builds on the renormalization scheme of Kozma and
Nitzan~\cite[\S4]{KozmaNitzan2024} and the seed-and-exploration method of Grimmett and
Marstrand~\cite{GrimmettMarstrand1990}. Figure~\ref{fig:roadmap} shows the structure. Assume, for
a contradiction, that $\theta(\beta)>0$ at $\beta=\beta_c$. The argument has four steps.
\begin{enumerate}
\item \emph{Faces at $\beta$.} A large box $\Lam_m(u)$ connects with high probability, inside a
surrounding cube $\Lam_\ell(u)$, to each of the $d2^d$ quarter-faces of that cube (a face cut in half
in every transverse coordinate). Only three radii $\ell$ are needed, so these finitely many estimates
persist at some $q<\beta$.
\item \emph{Extension.} If a vertex $o$ outside a box $B$ reaches $B$, then it reaches, with
almost the same probability, any target set $T$ that every point near $B$ can hit through a
quarter-face. The cluster of $o$ enters a neighbourhood of $B$ at many possible points; at each one a
small \emph{seed} box, if all its edges are open, gives a \emph{relay point} that reaches $T$ with high
probability. The gluing inequality (Theorem~\ref{thm:gluing}) turns ``every relay point reaches $T$ with
high probability'' into ``$o$ reaches $T$ if it reaches a relay point''.
\item \emph{Relay.} A fixed sequence of extensions moves a connection from the centre of one
coarse cube to the centre of a neighbouring one.
\item \emph{Exploration.} Cubes indexed by $\Z^2$ are explored one at a time, keeping a high
conditional connection probability for each pending cube; a planar contour count then yields an
infinite cluster at $q<\beta_c$, a contradiction.
\end{enumerate}

Long edges break two nearest-neighbour features of this argument.

\begin{figure}[tp]
\centering
\begin{tikzpicture}[scale=0.95, >={Stealth[length=1.6mm]}]
  % shell
  \fill[black!10] (-2.5,-2.5) rectangle (2.5,2.5);
  \fill[white] (-2,-2) rectangle (2,2);
  % highlighted outer thick quarter-face
  \fill[black!35] (2,0) rectangle (2.5,2.5);
  \draw[line width=0.6pt] (-2,-2) rectangle (2,2);
  \draw[black!50, dashed] (-2.5,-2.5) rectangle (2.5,2.5);
  % seed box
  \draw[line width=0.6pt, fill=black!20] (-0.35,-0.35) rectangle (0.35,0.35);
  \node at (0,-0.6) {\small $\Lambda_m(u)$};
  % open path with a short last edge
  \draw[line width=0.8pt, decorate, decoration={random steps, segment length=3.2pt, amplitude=1.1pt}]
    (0.25,0.2) -- (0.8,0.65) -- (1.3,0.45) -- (1.62,1.15) -- (1.86,1.3);
  \draw[line width=0.8pt] (1.86,1.3) -- (2.28,1.48);
  \fill (1.86,1.3) circle (1.3pt); \fill (2.28,1.48) circle (1.3pt);
  \node[anchor=west] at (2.6,1.48) {\small short last edge, length $\le w$};
  \node[anchor=west] at (2.6,0.55) {\small $F$, an outer thick quarter-face};
  \draw[->, black!70] (2.58,0.55) -- (2.3,0.8);
  % long edge
  \fill (-1.1,-0.75) circle (1.3pt); \fill (-3.3,-2.95) circle (1.3pt);
  \draw[line width=0.7pt, dashed] (-1.1,-0.75) -- (-3.3,-2.95);
  \node[anchor=north, align=center] at (-3.1,-3.05) {\small long edge: total probability\\ \small at most $\beta(2\ell+1)^dT(w)$};
  % labels
  \node[anchor=north west] at (-1.95,1.95) {\small $\Lambda_\ell(u)$};
  \draw[decorate, decoration={brace, amplitude=3pt}] (2,-2.6) -- (2.5,-2.6) node[midway, below=2pt] {};
  \node at (2.25,-2.95) {\small $w$};
  \node[anchor=south west, black!60] at (2.5,2.5) {\small $\Lambda_{\ell+w}(u)$};
\end{tikzpicture}
\caption{Lemma~\ref{lem:faces} in a planar cross-section. An open path from $\Lambda_m(u)$ that
leaves $\Lambda_\ell(u)$ by an edge of length at most $w$ lands in the shell
$\Lambda_{\ell+w}(u)\setminus\Lambda_\ell(u)$ (light grey), which the $N_F$ outer thick quarter-faces
cover; one is shown in dark grey. A path can also leave by a long edge (dashed), which a
nearest-neighbour argument would not see; such edges are rare when $\ell^dT(w)$ is small.}
\label{fig:faces}
\end{figure}
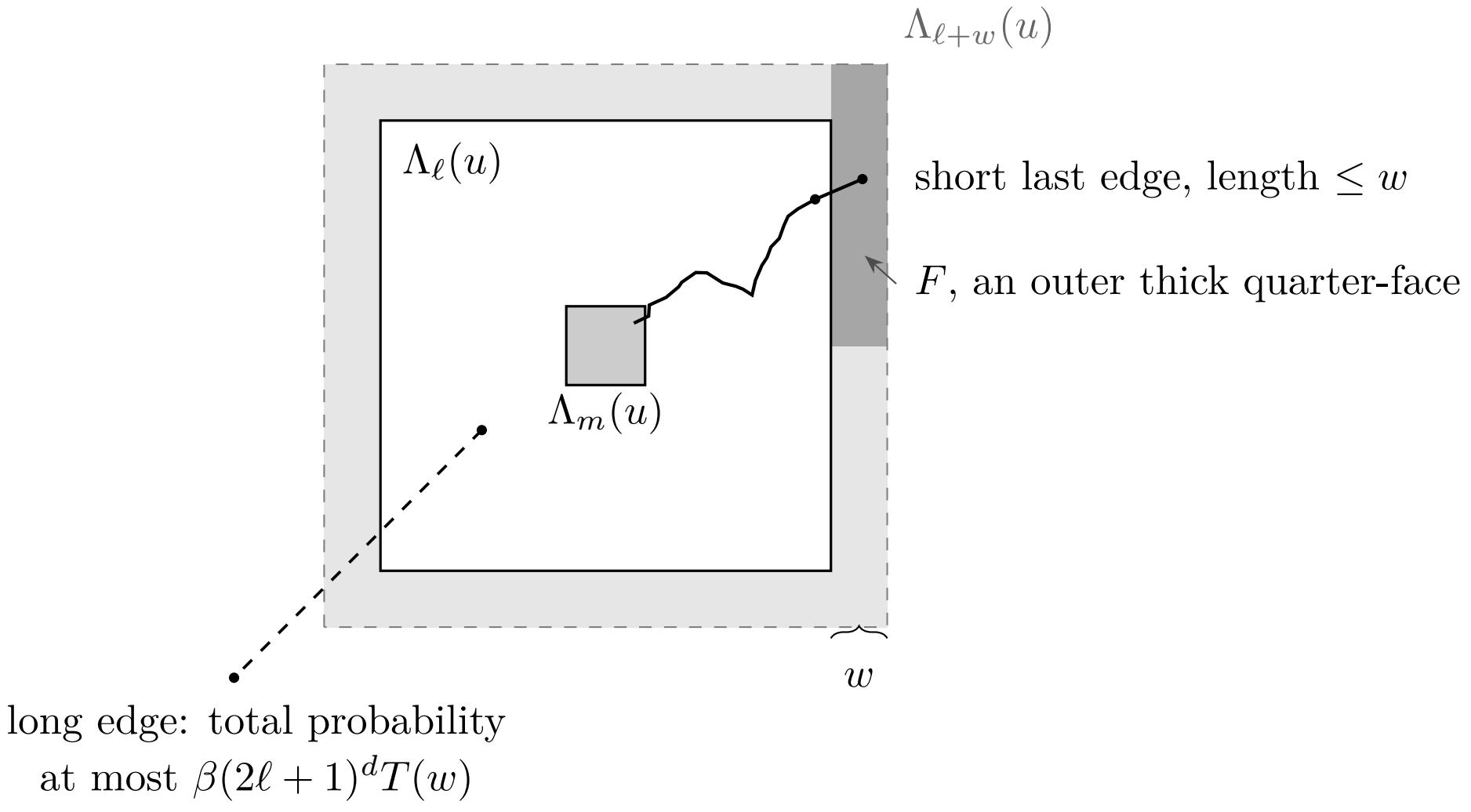

First, a path that leaves a box need not pass through its boundary: its first step outside can
be a long jump. In Step~1 we therefore use \emph{outer thick quarter-faces}, slabs of width $w$
just outside a face of the cube (Figure~\ref{fig:faces}), and we require the last edge of the
connection to have length at most $w$. Jumps longer than $w$ out of $\Lam_\ell$ have total
probability at most $\beta(2\ell+1)^dT(w)$. The tail condition makes this small for a $w$ that is
a small fixed fraction of $\ell$, so the thick faces still fit inside the relay geometry of Step~3.

\begin{figure}[tp]
\centering
\begin{tikzpicture}[scale=0.9, >={Stealth[length=1.6mm]}]
  % B^(R) and D
  \draw[black!40, dotted, line width=0.7pt] (-4.6,-3.3) rectangle (4.6,3.3);
  \node[anchor=south east, black!55] at (4.6,-3.3) {\small $B^{(R)}$};
  % shell S_j between C and B^(j+w)
  \fill[black!10] (-3.4,-2.5) rectangle (3.4,2.5);
  \fill[white] (-2.9,-2.0) rectangle (2.9,2.0);
  \draw[line width=0.7pt] (-2.9,-2.0) rectangle (2.9,2.0);
  \draw[black!50, dashed] (-3.4,-2.5) rectangle (3.4,2.5);
  \node[anchor=north west] at (-2.9,2.0) {\small $C=B^{(j)}$};
  \node[anchor=south east, black!60] at (3.4,2.5) {\small $B^{(j+w)}$};
  % B
  \draw[line width=0.6pt, fill=black!25] (-0.8,-0.5) rectangle (0.8,0.5);
  \node at (0,0) {\small $B$};
  % outside cluster: path from o to shell, with branches
  \fill (-6.2,0.4) circle (1.6pt); \node[anchor=east] at (-6.3,0.4) {\small $o$};
  \draw[line width=0.9pt, decorate, decoration={random steps, segment length=3.5pt, amplitude=1.2pt}]
    (-6.2,0.4) -- (-5.4,0.9) -- (-4.9,0.3) -- (-4.3,0.8) -- (-3.7,0.5) -- (-3.15,0.75);
  \draw[line width=0.9pt, decorate, decoration={random steps, segment length=3.5pt, amplitude=1.2pt}]
    (-4.9,0.3) -- (-4.5,-0.6) -- (-3.6,-1.0) -- (-3.1,-1.4);
  \draw[line width=0.9pt, decorate, decoration={random steps, segment length=3.5pt, amplitude=1.2pt}]
    (-5.4,0.9) -- (-5.6,2.2) -- (-4.9,2.9);
  \node[anchor=south] at (-5.6,2.3) {\small $\mathcal C$};
  % crossing edges and landing points with seed boxes
  \foreach \a/\b/\c/\e in {-3.15/0.75/-2.55/0.95, -3.1/-1.4/-2.3/-1.1, -3.7/0.5/-1.75/0.25}{
     \draw[line width=0.5pt] (\a,\b) -- (\c,\e);
     \fill (\c,\e) circle (1.2pt);
  }
  \draw[line width=0.5pt, fill=black!15, fill opacity=0.6] (-2.85,0.6) rectangle (-2.15,1.3);
  \draw[line width=0.5pt, fill=black!15, fill opacity=0.6] (-2.6,-1.45) rectangle (-1.9,-0.75);
  \draw[line width=0.5pt, fill=black!15, fill opacity=0.6] (-2.1,-0.1) rectangle (-1.4,0.6);
  % long jump over the shell
  \fill (-4.9,2.9) circle (1.2pt);
  \draw[line width=0.7pt, dashed] (-4.9,2.9) -- (0.4,1.45);
  \fill (0.4,1.45) circle (1.2pt);
  % annotations
  \node[anchor=west, align=left] at (4.85,1.2) {\small crossing edges from $\mathcal C$\\ \small into $C$; total weight $W_j$};
  \node[anchor=west, align=left] at (4.85,0.0) {\small seed boxes $U(z)$ at\\ \small landing points $z$};
  \node[anchor=west, align=left] at (4.85,-1.2) {\small shell $S_j$ of width $w$};
  \node[anchor=west, align=left] at (4.85,2.45) {\small dashed: jump over the shell,\\ \small controlled by $T(w)$};
  \draw[->, black!60] (4.8,-1.2) -- (3.15,-2.25);
  \draw[->, black!60] (4.8,0.0) -- (-1.4,0.45);
\end{tikzpicture}
\caption{The extension lemma (Lemma~\ref{lem:extension}) in a planar cross-section. The cluster
$\mathcal C$ of $o$ outside $C=B^{(j)}$ reaches the shell $S_j$. All crossing edges from
$\mathcal C$ into $C$ are closed with probability exactly $e^{-W_j}$, which is cheap when $W_j$ is
small; when $W_j$ is large, the landing points $z$ in one residue class carry disjoint seed boxes
$U(z)$, which give many independent chances to reach a good relay centre. Entrances that jump over
the shell are rare by the tail condition.}
\label{fig:extension}
\end{figure}
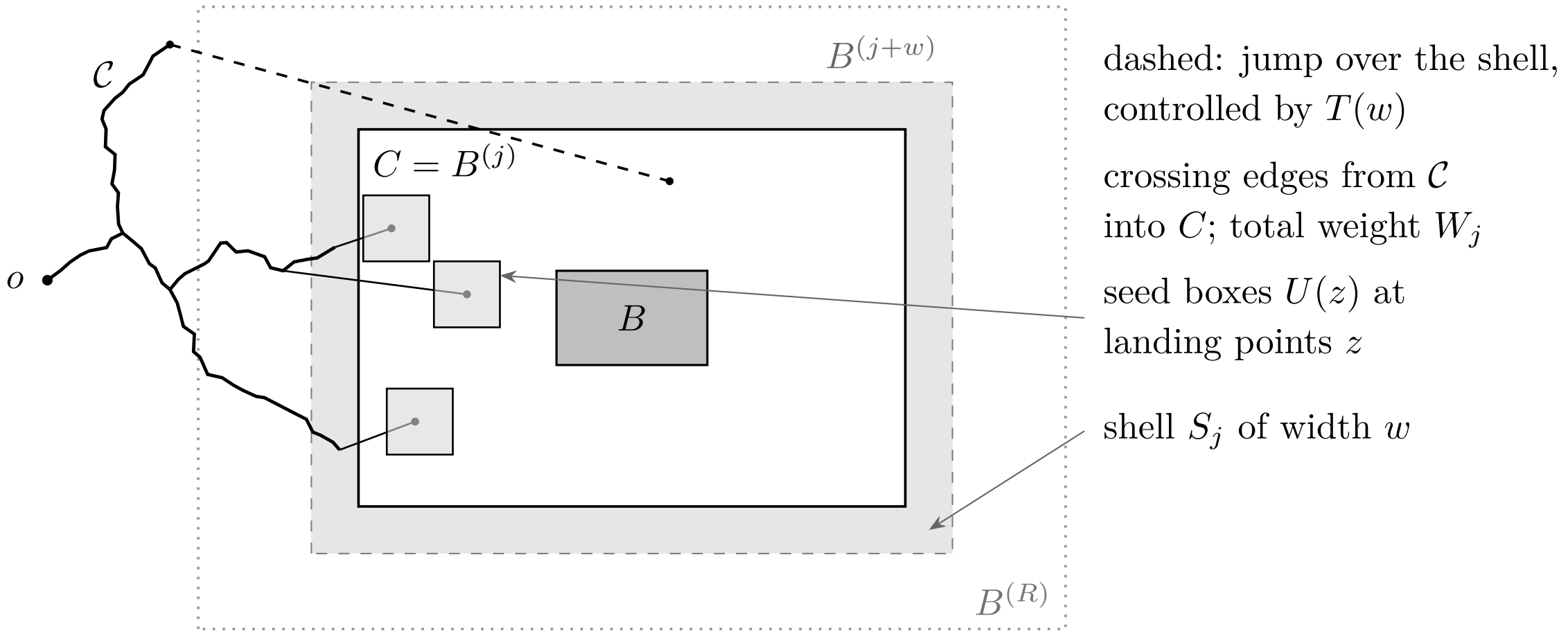

Second, in Step~2 the nearest-neighbour proof selects an expansion $B'=B^{(j)}$ of $B$ such that the
cluster of $o$ outside $B'$ has many exterior neighbours, using the fact that closing the single edge
from each exterior neighbour into $B'$ disconnects $o$ from $B'$. With long edges, every vertex of the
outside cluster has edges into $B'$. We replace the count of exterior neighbours by the
\emph{crossing weight}
\[
W=\sum_{x\in\cC}\ \sum_{z\in B'}\beta J(x-z),
\]
where $\cC$ is the cluster of $o$ in the complement of $B'$; the probability that all crossing edges
are closed is exactly $e^{-W}$ (Figure~\ref{fig:extension}). If $W$ is small this disconnection
is cheap, and disjointness of the corresponding events over well-separated layers shows that a
small weight is rare at some layer. If $W$ is large, the landing points of the crossing edges
support many independent attempts to open a seed box, because a single seed box $\Lam_m(v)$ can
absorb at most $\beta(2m+1)^d\norm{J}_1$ of the weight. Entrances that jump over a shell of width $w$
around $B'$ are again controlled by the tail condition. The other steps of~\cite{OpenAIZ3} need only
bookkeeping: $d$ contraction steps instead of three in the relay between neighbouring cubes, margins
enlarged by $w$, and an exploration whose root cube is opened without fixing the edges that leave it.

The scales are chosen in the order $m\le w\ll R\ll r$, where $R$ and $r$ are large constant
multiples of $w$ and the tail condition makes $\beta r^dT(w)$ small (Lemma~\ref{lem:scales}). Only
one such choice is needed; this is why a $\liminf$ suffices.

\emph{Adapted results.} Lemma~\ref{lem:faces} and Lemma~\ref{lem:scales} adapt~\cite[Lemma~3.3]{OpenAIZ3} (outer thick faces
and short exits); Lemma~\ref{lem:extension} adapts~\cite[Lemma~4.1]{OpenAIZ3}, itself adapted
from~\cite[Lemma~10]{KozmaNitzan2024} (crossing weight, layers spaced by $w$, seeds at landing
points); Lemma~\ref{lem:relay} adapts~\cite[Lemma~5.1]{OpenAIZ3} to $\Z^d$ and to the margins
$R+w$; Section~\ref{sec:exploration} adapts~\cite[\S6, Lemmas~6.1--6.3]{OpenAIZ3} (long-range edges between
neighbouring cubes, a hypothesis on edges entering the union of two cubes, and the root). The planar contour count is used as in
\cite[Lemma~3.1]{OpenAIZ3}, and we include its proof.

\subsection{Organization}
Section~\ref{sec:prelim} collects the inputs. Sections~\ref{sec:faces}--\ref{sec:exploration} carry
out Steps~1--4, and Section~\ref{sec:corollaries} proves the corollaries and records where each
hypothesis is used.

\section{Preliminaries}\label{sec:prelim}

\paragraph{Notation.} $\Lam_b(u)=u+([-b,b]^d\cap\Z^d)$ and $\Lam_b=\Lam_b(0)$. A \emph{rectangle}
is a nonempty product of integer intervals, and $B^{(j)}$ denotes the expansion of a rectangle $B$
by $j$ in every coordinate direction. Throughout, $\norm{J}_1=\sum_xJ(x)$; it is finite under (Tail),
since $T(n_0)<\infty$ for some $n_0$ and the remaining sum is finite.

\paragraph{Product laws on finite graphs.} In Sections~\ref{sec:extension}--\ref{sec:exploration}
we work on finite graphs $G=(V,E)$ with $V\subset\Z^d$, whose edges are pairs $\{x,y\}$ with
$J(x-y)>0$, with independent bits of individual probabilities $p_e\in[0,1]$. We put
$\lambda_e=-\log(1-p_e)\in[0,\infty]$, so that $\Pp(e\text{ open})\le\lambda_e$ and the probability
that all edges of a set $S$ are closed is $\exp(-\sum_{e\in S}\lambda_e)$. A bit is \emph{internal}
to a vertex set $U$ if both endpoints lie in $U$. For $q>0$ we say that $(F_q(U))$ holds if every
pair $\{x,y\}\subset U$ with $J(x-y)>0$ is an edge of $G$ with $\lambda_{xy}=qJ(x-y)$.

\paragraph{The gluing inequality.}
\begin{theorem}[\cite{OpenAIZ3,BlancRenaudie2026,Gouere2026}]\label{thm:gluing}
Consider independent bond percolation on a finite graph with arbitrary edge probabilities in
$[0,1]$. For vertices $s,t$ and a nonempty vertex set $A$,
\[
\Pp(s\lr A,\ s\lr t)\ \ge\ \Pp(s\lr A)\,\min_{a\in A}\Pp(a\lr t).
\]
\end{theorem}
This is a joint form of Conjecture~1 of Kozma and Nitzan~\cite{KozmaNitzan2024}, in which the event
$\{s\lr t\}$ alone appears on the left. Theorem~\ref{thm:gluing} is proved for hypergraphs, with $t$
replaced by a vertex set, in~\cite[Proposition~2.1]{OpenAIZ3}, and in the form above
in~\cite[Theorem~1]{Gouere2026}, where the extension from $(0,1)$ to $[0,1]$ is by continuity; it also
follows from~\cite[Theorem~1]{BlancRenaudie2026} together with Harris's inequality. The Lean
formalization of Leder~\cite{Leder2026} proves the additive form: if $\Pp(a\lr t)\ge1-b$ for every
$a\in A$, then $\Pp(s\lr t)\ge\Pp(s\lr A)-b$~\cite[Theorem~2.2]{Leder2026}. Either form gives the following
corollary, which is~\cite[Corollary~2.3]{OpenAIZ3}; we include the short proof from
Theorem~\ref{thm:gluing}.

\begin{corollary}[failure comparison]\label{cor:failure}
Consider independent bond percolation on a finite graph $G$ with arbitrary edge probabilities in
$[0,1]$, let $A,T$ be vertex sets and $b\in[0,1]$. If $\Pp(a\lr T)\ge1-b$ for every $a\in A$, then
$\Pp(o\lr A,\ o\nlr T)\le b$ for every vertex $o$. The same holds after any deterministic set of bits
has been fixed.
\end{corollary}
\begin{proof}
If $o\in T$, $A=\varnothing$ or $T=\varnothing$ (which forces $b=1$ when $A\ne\varnothing$) there is
nothing to prove. Otherwise let $G'$ be obtained from $G$ by
identifying the vertices of $T$ to a single vertex $t$, removing loops, and merging the parallel edges
$e_1,\dots,e_k$ between two vertices into one edge, open if and only if one of the $e_i$ is open; this
couples the two models, and the merged edge has probability $1-\prod_i(1-p_{e_i})$. Under the coupling,
$o\nlr T$ in $G$ if and only if $o\nlr t$ in $G'$, and then the open clusters of $o$ in $G$ and $G'$
coincide. Let $A'$ be the image of $A$ in $G'$. Hence
$\Pp_G(o\lr A,\ o\nlr T)=\Pp_{G'}(o\lr A',\ o\nlr t)$. In $G'$ every $a\in A'$ satisfies
$\Pp_{G'}(a\lr t)\ge1-b$: for $a\notin T$ this is $\Pp_G(a\lr T)\ge1-b$, and $t\lr t$ always. Theorem~\ref{thm:gluing}
in $G'$ gives
$\Pp_{G'}(o\lr A',\ o\nlr t)=\Pp_{G'}(o\lr A')-\Pp_{G'}(o\lr A',\ o\lr t)\le b\,\Pp_{G'}(o\lr A')\le b$.
Fixing a bit means giving it probability $0$ or $1$.
\end{proof}

We also use Harris's inequality~\cite{Harris1960}: increasing events of a finite product of
independent bits are positively correlated, and so are decreasing events.

\paragraph{A planar boundary count.} The following lemma is~\cite[Lemma~3.1]{OpenAIZ3}; we include
the short proof.
\begin{lemma}\label{lem:contours}
There are families $\mathcal C_n$ of simple cycles of length $n$ in the dual lattice
$\Z^2+(\frac12,\frac12)$, with $|\mathcal C_n|\le2n3^n$, such that the dual edge boundary of every
finite $X\subset\Z^2$ containing $0$ contains a cycle in $\bigcup_n\mathcal C_n$. For each cycle in
$\mathcal C_n$ one can fix at least $n/7$ of its crossed primal edges with pairwise disjoint
endpoints.
\end{lemma}
\begin{proof}
Each dual vertex has even degree in the dual boundary of $X$ (it counts membership changes around
a primal square), so the boundary decomposes into edge-disjoint simple cycles. Along the positive
horizontal ray from $0$ the number of boundary crossings is odd, so some cycle of the decomposition
crosses it an odd number of times; a closed cycle crosses the whole horizontal line an even number
of times, so this cycle also crosses the negative ray. If its length is $n$, some positive crossing
has horizontal coordinate $t+\frac12$ with $0\le t<n$. Let $\mathcal C_n$ be the set of simple
cycles of length $n$ with such an anchored crossing and an odd number of positive crossings: there
are at most $n$ anchors, two orientations and $3^n$ non-backtracking continuations. Finally select
crossed primal edges greedily, discarding those incident to an endpoint of a selected edge; each
selection discards at most seven edges, since primal degrees are four.
\end{proof}

\paragraph{Elementary facts.}
\begin{lemma}\label{lem:basic}
Under \textup{(Sym)}, \textup{(Pos)} and $\norm{J}_1<\infty$, if $d\ge2$ then $0<\beta_c<\infty$.
If $\theta(\beta)>0$, then an infinite cluster exists almost surely, and
\begin{equation}\label{eq:seed-infinite}
\Pp_\beta(\Lam_m\text{ meets an infinite cluster})\to1\qquad(m\to\infty).
\end{equation}
\end{lemma}
\begin{proof}
The nearest-neighbour pairs in a coordinate plane are open independently with probability
$1-e^{-\beta J(e_1)}$, by (Sym) and (Pos), and planar bond percolation with parameter close to $1$
percolates by Lemma~\ref{lem:contours} and a union bound over dual cycles; so $\beta_c<\infty$. The
cluster of $0$ is dominated by a Galton--Watson tree with mean offspring
$\sum_y(1-e^{-\beta J(y)})\le\beta\norm{J}_1$, so $\beta_c\ge1/\norm{J}_1$. Removing an edge
splits a connected component into at most two components, so the event that an infinite cluster
exists is unchanged by modifying finitely many bits and is a tail event; by Kolmogorov's
$0$--$1$ law it has probability $1$ when $\theta(\beta)>0$. The events that $\Lam_m$ meets an
infinite cluster increase to this event.
\end{proof}

\section{Faces at a percolating parameter, and the choice of scales}\label{sec:faces}

\begin{definition}\label{def:faces}
Let $u\in\Z^d$ and $\ell,w\ge1$. For $i\in\{1,\dots,d\}$, $\sigma\in\{\pm1\}$ and
$\varepsilon=(\varepsilon_j)_{j\ne i}\in\{\pm1\}^{d-1}$, the \emph{outer thick quarter-face}
$F_{i,\sigma,\varepsilon}(u,\ell)$ is
\[
\bigl\{u+x:\ \sigma x_i\in(\ell,\ell+w],\ 0\le\varepsilon_jx_j\le\ell+w\ (j\ne i)\bigr\}.
\]
There are $N_F=d\,2^d$ of them. For $m<\ell$ and such a face $F$, let $X(u,\ell,F)$ be the event
that there is an open path $x_0,\dots,x_k$ with $x_0\in\Lam_m(u)$, $x_0,\dots,x_{k-1}\in\Lam_\ell(u)$,
$x_k\in F$ and $\norm{x_k-x_{k-1}}\le w$.
\end{definition}

The event $X(u,\ell,F)$ is increasing. It depends only on the bits internal to $\Lam_\ell(u)$ and
on the pairs between $\Lam_\ell(u)$ and $F$, all of which are internal to $\Lam_{\ell+w}(u)$.

\begin{lemma}[faces at a percolating parameter]\label{lem:faces}
Assume \textup{(Sym)} and $\theta(\beta)>0$. For $\ell>m$, $w\ge1$, every $u$ and every outer thick
quarter-face $F$ of $(u,\ell)$,
\[
\Pp_\beta\bigl(X(u,\ell,F)^c\bigr)\le\Bigl[\Pp_\beta(\Lam_m\text{ meets no infinite cluster})
+\beta(2\ell+1)^dT(w)\Bigr]^{1/N_F}.
\]
\end{lemma}
\begin{proof}
Take $u=0$. Let $E_{\mathrm{long}}$ be the event that some open edge $\{x,y\}$ has $x\in\Lam_\ell$
and $\norm{x-y}>w$; then
$\Pp(E_{\mathrm{long}})\le\sum_{x\in\Lam_\ell}\sum_{\norm{z}>w}(1-e^{-\beta J(z)})\le\beta(2\ell+1)^dT(w)$.
If $\Lam_m$ meets an infinite cluster, take a finite open path from $\Lam_m$ to a vertex outside
$\Lam_\ell$, and let $y$ be its first vertex outside $\Lam_\ell$ and $x$ its predecessor. Off
$E_{\mathrm{long}}$ we have $\norm{x-y}\le w$, so $y\in\Lam_{\ell+w}\setminus\Lam_\ell$. Choose $i$ with
$|y_i|>\ell$, $\sigma=\operatorname{sign}y_i$ and $\varepsilon_j=\operatorname{sign}y_j$ (either sign
if $y_j=0$); then $y\in F_{i,\sigma,\varepsilon}(0,\ell)$ and the path up to $y$ witnesses
$X(0,\ell,F_{i,\sigma,\varepsilon})$. Hence
$\Pp(\bigcap_FX(0,\ell,F)^c)\le\Pp(\Lam_m\text{ meets no infinite cluster})+\Pp(E_{\mathrm{long}})$.
The events $X(0,\ell,F)^c$ are decreasing and depend on finitely many bits. The signed coordinate
permutations fixing $0$ act transitively on the faces and, by (Sym), preserve the law, $\Lam_m$
and $\Lam_\ell$; so all $N_F$ complements have the same probability $a$. Harris's inequality gives
$a^{N_F}\le\Pp(\bigcap_FX(0,\ell,F)^c)$.
\end{proof}

\begin{lemma}[choice of scales]\label{lem:scales}
Assume \textup{(Sym)}, \textup{(Tail)} and $\theta(\beta)>0$. There are $f,\delta,\eta,\tau,\gamma\in(0,1)$,
positive integers $m,K,w,R,r$, constants $\rho_0,A,c_A,L>0$ and $q\in(\beta/2,\beta)$ such that:
\begin{enumerate}
\item[(S1)] $f<\frac12$ and $\sum_{n\ge1}2n3^n(2f)^{n/7}<1$;
\item[(S2)] $\delta+4(d+10)\eta/\delta<f$, $\tau<\eta/4$, $\gamma<\delta$ and $2\gamma/\tau<\eta/4$;
\item[(S3)] $\Pp_\beta(\Lam_m\text{ meets no infinite cluster})<\frac12(\gamma/2)^{N_F}$;
\item[(S4)] $\rho_0=\prod_e(1-e^{-(\beta/2)J(e)})$, the product over the pairs $e$ internal to $\Lam_m$ with $J(e)>0$;
$A=\beta\norm{J}_1$; $c_A=(1-e^{-A})/A$; $4/(K\rho_0)<\eta/4$; and $L=K(4m+3)^d/c_A$;
\item[(S5)] $\beta(41r+1)^dT(w)\le\min\{\frac12(\gamma/2)^{N_F},\eta/8\}$;
\item[(S6)] $w\ge m$, $R\ge m+w+8e^Lw/\eta$, $r>4(d+12)(R+w+1)$ and $r>m+w$;
\item[(S7)] for $\ell\in\cL:=\{2r,2r+1,10r-w\}$, every $u$ and every outer thick quarter-face $F$ of
$(u,\ell)$, $\Pp_q(X(u,\ell,F))>1-\gamma$.
\end{enumerate}
\end{lemma}
\begin{proof}
Choose $f$, then $\delta,\eta,\tau,\gamma$ in this order; then $m$ by~\eqref{eq:seed-infinite};
then $\rho_0,A,c_A,K,L$. None of these depends on $w$ or $r$. By (Tail) there are arbitrarily large
$w\ge m$ with $w^dT(w)$ as small as we like. For such $w$ put
$R=m+w+\lceil8e^Lw/\eta\rceil$ and $r=4(d+12)(R+w+2)$; then (S6) holds and $r\le C_0w$ with
$C_0$ depending only on $d,m,L,\eta$, so $\beta(41r+1)^dT(w)\le\beta(42C_0)^dw^dT(w)$ and (S5) holds
for a suitable such $w$. Every $\ell\in\cL$ satisfies $\ell>m$ and $2\ell+1\le41r+1$, so
Lemma~\ref{lem:faces}, (S3) and (S5) give $\Pp_\beta(X(u,\ell,F)^c)<\gamma/2$. These finitely many
events (up to translation) depend on finitely many bits whose probabilities are continuous in the
parameter, so (S7) holds for some $q\in(\beta/2,\beta)$.
\end{proof}

We fix these choices from now on, and write $\rho=\prod_e(1-e^{-qJ(e)})\ge\rho_0$ (over the same pairs) for the probability
that all pairs $e$ internal to a translate of $\Lam_m$ with $J(e)>0$ are open.

\section{Extension through a fresh region}\label{sec:extension}

\begin{lemma}[extension]\label{lem:extension}
Let $G=(V,E)$ be a finite graph as in Section~\ref{sec:prelim}, $D\subset V$ a rectangle with
$|D|\le(41r+1)^d$ and $(F_q(D))$, $B$ a rectangle with $B^{(R)}\subset D$, $o\in V\setminus B^{(R)}$
and $T\subset V$. Assume:
\begin{enumerate}
\item[(E2)] every edge $\{x,z\}$ of $G$ with $z\in B^{(R)}$ and $x\notin D$ has $\lambda_{xz}\le qJ(x-z)$;
\item[(E3)] for every $u\in B^{(R)}$ there are $\ell\in\cL$ and an outer thick quarter-face $F$ of $(u,\ell)$
with $\Lam_{\ell+w}(u)\subset D$ and $F\subset T$.
\end{enumerate}
Then $\Pp(o\lr B,\ o\nlr T)<\eta$.
\end{lemma}

\begin{proof}
By $(F_q(D))$ and (E2), every edge $\{x,z\}$ with $z\in B^{(R)}$ and $x\notin B^{(R)}$ has
$\lambda_{xz}\le qJ(x-z)$; hence $\sum_x\lambda_{xz}\le q\norm{J}_1\le A$ for $z\in B^{(R)}$.

\emph{Step 1: a layer.} For an integer $j\in[m,R)$ let $\zeta_j$ be the bits internal to
$V\setminus B^{(j)}$, let $\cC_j$ be the open cluster of $o$ in the graph on $V\setminus B^{(j)}$ with
these bits, and set
\[
W_j=\sum_{x\in\cC_j}\sum_{z\in B^{(j)}}\lambda_{xz},\qquad S_j=B^{(j+w)}\setminus B^{(j)},
\]
\[
L_j=\{W_j<L,\ \cC_j\cap S_j\ne\varnothing\},\qquad
\mathrm{Cl}_j=\{\text{every edge between }\cC_j\text{ and }B^{(j)}\text{ is closed}\}.
\]
Here $W_j\le|B^{(j)}|A<\infty$. The crossing bits between $V\setminus B^{(j)}$ and $B^{(j)}$ are
independent of $\zeta_j$, so $\Pp(\mathrm{Cl}_j\mid\zeta_j)=e^{-W_j}\ge e^{-L}$ on $L_j$. On
$\mathrm{Cl}_j$ the open cluster $C(o)$ of $o$ in $G$ equals $\cC_j$, because every other edge
leaving $\cC_j$ is internal to $V\setminus B^{(j)}$ and therefore closed by maximality.

Let $\cJ=\{m+kw:k\ge0\}\cap[m,R)$. If $j'<j$ both lie in $\cJ$, then $j'+w\le j$; on
$L_j\cap\mathrm{Cl}_j$ we have $C(o)=\cC_j\subset V\setminus B^{(j)}$, hence $\cC_{j'}=C(o)$, while
$S_{j'}\subset B^{(j)}$, so $\cC_{j'}\cap S_{j'}=\varnothing$ and $L_{j'}$ fails. The events
$L_j\cap\mathrm{Cl}_j$, $j\in\cJ$, are therefore pairwise disjoint, and
\[
e^{-L}\sum_{j\in\cJ}\Pp(L_j)\le\sum_{j\in\cJ}\Pp(L_j\cap\mathrm{Cl}_j)\le1.
\]
Since $|\cJ|\ge(R-m)/w$, we can fix $j\in\cJ$ with $\Pp(L_j)\le e^Lw/(R-m)\le\eta/8$, by (S6).

Let $E_{\mathrm{long}}$ be the event that some open edge $\{x,z\}$ has $z\in B^{(R)}$ and
$\norm{x-z}>w$. Then
$\Pp(E_{\mathrm{long}})\le\sum_{z\in B^{(R)}}\sum_{\norm{x-z}>w}\lambda_{xz}\le q|D|T(w)\le\eta/8$
by (S5). If $o\lr B$, $W_j<L$ and $\cC_j\cap S_j=\varnothing$, then an open path from $o$ enters
$B^{(j)}$ for the first time at some $z$ from a predecessor $x\in\cC_j\setminus B^{(j+w)}$, so
$\norm{x-z}>w$ and $E_{\mathrm{long}}$ occurs. Hence
\begin{equation}\label{eq:low-weight}
\Pp(o\lr B,\ W_j<L)\le\Pp(L_j)+\Pp(E_{\mathrm{long}})\le\eta/4.
\end{equation}

\emph{Step 2: relay quality.} Put $C=B^{(j)}$ and let $\xi$ be the bits internal to $C$. Every side
of $C$ has at least $2j+1\ge2m+1$ points. For $z\in C$ let $v(z)$ be obtained from $z$ by clamping
each coordinate to the corresponding side of $C$ shrunk by $m$ at both ends, let
$U(z)=\Lam_m(v(z))\subset C$ (so $z\in U(z)$), let $O(z)$ be the event that all bits internal to
$U(z)$ are open, so $\Pp(O(z))=\rho$, and put
\[
g(z)=\Pp(v(z)\lr T\mid\xi),\qquad A_\xi=\{v(z):z\in C,\ g(z)>1-\tau\}.
\]
$A_\xi$ is a function of $\xi$. Fix $z$. By (E3) at $u=v(z)\in B^{(R)}$ there are $\ell\in\cL$ and $F$
with $\Lam_{\ell+w}(v(z))\subset D$ and $F\subset T$; the event $X(v(z),\ell,F)$ uses only bits
internal to $D$, so by $(F_q(D))$ and (S7) its probability exceeds $1-\gamma$. By (Pos) and (Sym)
an open box is connected through its nearest-neighbour pairs, so on $O(z)\cap X(v(z),\ell,F)$ we
have $v(z)\lr T$. Harris's inequality gives $\Pp(O(z),\ v(z)\nlr T)\le\rho\gamma$, and since $O(z)$
is $\xi$-measurable,
\begin{equation}\label{eq:relay-quality}
\tau\,\Pp\bigl(O(z),\ g(z)\le1-\tau\bigr)\le\E\bigl[\1_{O(z)}(1-g(z))\bigr]=\Pp\bigl(O(z),\ v(z)\nlr T\bigr)\le\rho\gamma.
\end{equation}

\emph{Step 3: entrances.} Condition on $\zeta=\zeta_j$; write $\cC=\cC_j$ and
$a_z=\sum_{x\in\cC}\lambda_{xz}$ for $z\in C$, so that $W_j=\sum_{z\in C}a_z$ and $a_z\le A$. Among the
$(4m+3)^d$ residue classes of $C$ modulo $4m+3$, choose (as a function of $\zeta$) a class $Z$ with
$\sum_{z\in Z}a_z\ge W_j/(4m+3)^d$. Distinct $z,z'\in Z$ satisfy $\norm{z-z'}\ge4m+3$, and clamping
moves each coordinate by at most $m$, so $\norm{v(z)-v(z')}\ge2m+3$ and the boxes $U(z)$, $z\in Z$,
are disjoint.

Let $\Sigma_z$ be the event that some edge between $\cC$ and $z$ is open and $O(z)$ holds. Given
$\zeta$, the crossing bits are independent of $\xi$, edges into distinct $z$ are distinct, and the
events $O(z)$, $z\in Z$, use disjoint bits of $\xi$; hence the events $\Sigma_z$, $z\in Z$, are
independent with $\Pp(\Sigma_z\mid\zeta)=\rho(1-e^{-a_z})$. Let $S=\sum_{z\in Z}\1_{\Sigma_z}$ and
$\mu=\sum_{z\in Z}(1-e^{-a_z})$. Since $a\mapsto(1-e^{-a})/a$ is decreasing, $1-e^{-a_z}\ge c_Aa_z$
and $\mu\ge c_AW_j/(4m+3)^d$, which is at least $K$ on $\{W_j\ge L\}$. By Chebyshev's inequality,
on $\{W_j\ge L\}$,
\[
\Pp(S<\rho\mu/2\mid\zeta)\le\frac{4}{\rho\mu}\le\frac4{\rho_0K}<\frac\eta4 .
\]
Let $Z_{\mathrm{bad}}=\sum_{z\in Z}\1_{\Sigma_z}\1\{g(z)\le1-\tau\}$. By the independence of the
crossing bits from $\xi$ and~\eqref{eq:relay-quality},
$\E[Z_{\mathrm{bad}}\mid\zeta]=\sum_z(1-e^{-a_z})\Pp(O(z),g(z)\le1-\tau)\le\mu\rho\gamma/\tau$, so
$\Pp(Z_{\mathrm{bad}}\ge\rho\mu/2\mid\zeta)\le2\gamma/\tau<\eta/4$. If $S\ge\rho\mu/2>Z_{\mathrm{bad}}$,
some $z\in Z$ has $\Sigma_z$ and $g(z)>1-\tau$: then $o$ reaches some $x\in\cC$, the open edge
$\{x,z\}$, and $v(z)\in A_\xi$ through the open box $U(z)$. Together with~\eqref{eq:low-weight},
\[
\Pp(o\lr B,\ o\nlr A_\xi)\le\tfrac34\eta .
\]

\emph{Step 4: the target.} Condition on $\xi$. Then $A_\xi$ is fixed, the remaining bits keep their
product law, and $\Pp(a\lr T\mid\xi)>1-\tau$ for every $a\in A_\xi$. Corollary~\ref{cor:failure},
applied with the bits of $\xi$ fixed, gives $\Pp(o\lr A_\xi,\ o\nlr T\mid\xi)\le\tau<\eta/4$.
Splitting according to whether $o\lr A_\xi$ gives $\Pp(o\lr B,\ o\nlr T)<\eta$.
\end{proof}

\section{Relay between neighbouring coarse cubes}\label{sec:relay}

For $v=(v_1,v_2)\in\Z^2$ let
\[
c_v=(20r+1)(v_1,v_2,0,\dots,0)\in\Z^d,\qquad Q_v=\Lam_{10r}(c_v),\qquad B_v=\Lam_{4r}(c_v).
\]
The cubes $Q_v$ are vertex-disjoint. For neighbours $v,x$ in $\Z^2$, let $I_{vx}$ be the finite set
of pairs $\{y,y'\}$ with $y\in Q_v$, $y'\in Q_x$ and $J(y-y')>0$. The union $D_{vx}=Q_v\cup Q_x$ is a
rectangle with $|D_{vx}|=(40r+2)(20r+1)^{d-1}\le(41r+1)^d$, whose internal bits are those internal to
$Q_v$, those internal to $Q_x$, and $I_{vx}$. Let $D^\circ_{vx}$ be the set of points of $D_{vx}$ at
distance more than $2r$ from $\Z^d\setminus D_{vx}$.

\begin{lemma}[relay]\label{lem:relay}
Let $v,x$ be neighbours in $\Z^2$ and $G$ a finite graph as in Section~\ref{sec:prelim} with
$D=D_{vx}\subset V$ and $(F_q(D))$. Assume
\begin{enumerate}
\item[(E2$^\circ$)] every edge $\{y,z\}$ of $G$ with $z\in D^\circ_{vx}$ and $y\notin D$ has
$\lambda_{yz}\le qJ(y-z)$.
\end{enumerate}
If $o\in V\setminus D$, then $\Pp(o\lr B_v,\ o\nlr B_x)\le(d+10)\eta$.
\end{lemma}

\begin{proof}
By (Sym) we may assume $c_v=0$ and $c_x=(20r+1)e_1$, so that
$D=[-10r,30r+1]\times[-10r,10r]^{d-1}$. Put $R'=R+w$. We define rectangles $P_0,\dots,P_{d+10}$ and
apply Lemma~\ref{lem:extension} with $B=P_{k-1}$ and $T=P_k$; the radius is $2r$ for $k\le d+9$
and $2r+1$ for $k=d+10$. Let $P_0=B_v$. For $t=1,\dots,d$ let $P_t$ be centred at $0$ with
half-width $2r+tR'$ in directions $h\le t$ and $4r+tR'$ in directions $h>t$. For $i=1,\dots,10$
put $s_i=2r$ ($i\le9$), $s_{10}=2r+1$, $z_i=(s_1+\dots+s_i)e_1$, $w_i=2r+(d+i)R'$, and
$P_{d+i}=\Lam_{w_i}(z_i)$.

\emph{Contraction steps.} Let $1\le t\le d$ and $u\in P_{t-1}^{(R)}$. Then $|u_h|\le2r+tR'-w$ for
$h<t$ and $|u_h|\le4r+tR'-w$ for $h\ge t$. Take $F=F_{t,\sigma,\varepsilon}(u,2r)$ with
$\sigma=-\operatorname{sign}u_t$ and $\varepsilon_j=-\operatorname{sign}u_j$ (with
$\operatorname{sign}0=+1$). A point $y\in F$ has $y_t=u_t+\sigma s$ with $s\in(2r,2r+w]$, so
$|y_t|\le|u_t|-2r\le2r+tR'$ if $|u_t|\ge2r+w$ and $|y_t|\le2r+w$ otherwise. For $j\ne t$, $y_j$ lies
between $u_j$ and $u_j-\operatorname{sign}(u_j)(2r+w)$, so $|y_j|\le\max\{|u_j|,2r+w\}$, which is at
most the half-width of $P_t$ in direction $j$. Thus $F\subset P_t$. The coordinates of
$\Lam_{2r+w}(u)$ are at most $6r+dR'+w<8r$ in absolute value, by (S6), so $\Lam_{2r+w}(u)\subset Q_v$
and $P_{t-1}^{(R)}\subset D^\circ_{vx}$. After the last contraction, $P_d=\Lam_{2r+dR'}(0)$.

\emph{Translation steps.} Let $1\le i\le10$ and $u\in P_{d+i-1}^{(R)}$, so
$\norm{u-z_{i-1}}\le w_{i-1}+R$ with $w_0=2r+dR'$. Take $F=F_{1,+,\varepsilon}(u,s_i)$ with
$\varepsilon_j=-\operatorname{sign}u_j$. For $y\in F$,
$y_1-(z_i)_1=(u_1-(z_{i-1})_1)+(y_1-u_1-s_i)$ with $y_1-u_1-s_i\in(0,w]$, so
$|y_1-(z_i)_1|\le w_{i-1}+R+w=w_i$, and $|y_j|\le\max\{|u_j|,s_i+w\}\le w_i$ for $j\ge2$. Thus
$F\subset P_{d+i}$. The first coordinates of $\Lam_{s_i+w}(u)$ lie in
$[-4r-(d+11)R'-1,\ 22r+(d+11)R'+1]$ and the others have absolute value at most $4r+(d+11)R'+1$; by
(S6) these lie in $D$ at distance more than $2r$ from its complement. Finally $z_{10}=c_x$ and
$w_{10}=2r+(d+10)R'<4r$, so $P_{d+10}\subset B_x$.

Each $P_{k-1}^{(R)}$ lies in $D^\circ_{vx}$, so (E2$^\circ$) gives (E2), and (E3) was just verified; also
$o\notin D$. If $o\lr P_0$ but $o\nlr P_{d+10}$, then $o\lr P_{k-1}$ and $o\nlr P_k$ for some $k$, and
the union bound over the $d+10$ applications of Lemma~\ref{lem:extension} proves the claim.
\end{proof}

\section{Exploration and the proof of Theorem~\ref{thm:main}}\label{sec:exploration}

We follow the exploration of~\cite[\S6]{OpenAIZ3}, with the long-range interfaces $I_{vx}$. Let
$o=0\in Q_0$, and let $\cI$ be the event that all bits internal to $Q_0$ are open. It prescribes
finitely many bits, each of probability in $(0,1)$, so $\Pp_q(\cI)>0$, and under
$\Pp_\cI=\Pp_q(\cdot\mid\cI)$ all other bits are independent with their parameter-$q$ law.

\paragraph{The exploration.} The coarse root $0\in\Z^2$ is declared good and visited, but is never queued
or processed; its four neighbours are put in a queue with parent $0$. A coarse vertex is \emph{visited} once it has been
queued, and \emph{pending} until it is processed; vertices are processed in queue order, and none is queued
twice. For a pending $v$, let $G(v)$ be the graph consisting of the cubes on the parent chain of
$v$ back to the root, the bits internal to these cubes, and the interfaces $I$ between successive
cubes of the chain, and no other pairs. Put $H_v=\{o\lr B_v\text{ in }G(v)\}$. To process $v$, reveal
the bits internal to $Q_v$ and those of $I_{pv}$, where $p$ is the parent of $v$. For each unvisited
neighbour $x$ of $v$, let $G'=G(v)\cup Q_x\cup I_{vx}$ (with the bits internal to $Q_x$) and compute
\[
b_x=\Pp_\cI(o\lr B_x\text{ in }G'\mid\text{history after this reveal}),
\]
which reveals nothing. Declare $v$ \emph{good} if $H_v$ holds and $b_x>1-\delta$ for every such $x$,
and in that case queue these $x$ with parent $v$; otherwise declare $v$ \emph{bad}. We maintain
\begin{equation}\label{eq:promise}
\Pp_\cI(H_v\mid\text{revealed history})>1-\delta\qquad\text{for every pending }v.
\end{equation}

\begin{lemma}[freshness and the initial promise]\label{lem:fresh}
Conditional on any revealed history of positive probability, the bits that are neither prescribed
by $\cI$ nor revealed are independent with their parameter-$q$ law, and~\eqref{eq:promise} holds at
every step.
\end{lemma}
\begin{proof}
Each step reveals a finite set of bits determined by earlier revealed values, and all decisions are
functions of revealed values; induction gives the first claim. The unknown bits of $H_x$ for a
pending $x$ with parent $v$ are internal to $Q_x$ or lie in $I_{vx}$. The cubes are vertex-disjoint,
each step reveals the internal bits of one cube and its own parent interface, and $G(x)$ contains no
other pairs, so these bits stay unrevealed while $x$ is pending; together with the definition of
\emph{good}, this keeps~\eqref{eq:promise} for every queued $x$ once it holds when $x$ is queued.

It remains to check~\eqref{eq:promise} for the root's neighbours. Let $v$ be one; by symmetry $Q_0$
lies in direction $-e_1$ from $Q_v$. Let $\ell'=10r-w\in\cL$, $u=c_v-we_1$ and
$F=F_{1,-,(+,\dots,+)}(u,\ell')$. Then $\Lam_{\ell'}(u)\subset Q_v$; $\Lam_m(u)\subset B_v$ since
$m+w\le4r$; and $F\subset Q_0$, since its first coordinates lie in $[(c_v)_1-10r-w,(c_v)_1-10r-1]$
and its other coordinates in $[0,10r]$. On $X(u,\ell',F)$, the box $\Lam_m(u)$ is joined inside $Q_v$
to a vertex with an open edge of length at most $w$ into $Q_0$, and this edge lies in $I_{0v}$, so
$X(u,\ell',F)\cap\cI\subset H_v$ (the cube $Q_0$ is open and connected by (Pos)). The event
$X(u,\ell',F)$ depends only on bits internal to $Q_v$ and on pairs of $I_{0v}$, which are independent
of $\cI$; hence $\Pp_\cI(H_v)\ge\Pp_q(X(u,\ell',F))>1-\gamma>1-\delta$ by (S7) and (S2).
\end{proof}

\begin{lemma}\label{lem:bad}
Immediately before a pending $v$ is processed, the conditional probability that $v$ is declared bad
is less than $f$.
\end{lemma}
\begin{proof}
Fix an unvisited neighbour $x$ of $v$ and let $D=D_{vx}$. Neither $v$ (pending) nor $x$ (unvisited)
has been processed and $I_{vx}$ is not a parent interface yet, so every bit internal to $D$ is fresh
with parameter $q$, and $(F_q(D))$ holds in $G'$. The pairs of $G'$ with one endpoint in $D$ and the
other outside are those of $I_{pv}$, which is still unrevealed, with $\lambda=qJ$; so (E2$^\circ$)
holds. Since $o\notin D$, Lemma~\ref{lem:relay} gives
$\Pp_\cI(H_v\cap\{o\nlr B_x\text{ in }G'\}\mid\text{history})\le(d+10)\eta$, and the tower property
gives $\Pp_\cI(H_v,\ b_x\le1-\delta\mid\text{history})\le(d+10)\eta/\delta$. With at most four such
$x$ and~\eqref{eq:promise}, $\Pp_\cI(v\text{ bad}\mid\text{history})\le\delta+4(d+10)\eta/\delta<f$
by (S2).
\end{proof}

\begin{lemma}\label{lem:badset}
For every deterministic set of $k$ distinct coarse vertices, the $\Pp_\cI$-probability that all of
them are processed and declared bad is at most $f^k$.
\end{lemma}
\begin{proof}
The root is good, so assume it is not in the set $S$. Let $\tau_j$ be the time at which the $j$th
member of $S$ is processed ($\tau_j=\infty$ if there is none) and $E_j$ the event that $\tau_j<\infty$
and the first $j$ processed members are bad. For each finite $t$ the event
$\{\tau_j=t\}\cap E_{j-1}$ is determined before the reveal at time $t$, so Lemma~\ref{lem:bad} gives
$\Pp_\cI(E_j)\le f\sum_t\Pp_\cI(\tau_j=t,E_{j-1})\le f\,\Pp_\cI(E_{j-1})$, and induction finishes the proof.
\end{proof}

\begin{proof}[Proof of Theorem~\ref{thm:main}]
By Lemma~\ref{lem:basic}, $0<\beta_c<\infty$. Suppose that $\theta(\beta_c)>0$ and choose the
parameters of Lemma~\ref{lem:scales} with $\beta=\beta_c$. Let $X\subset\Z^2$ be the set of good coarse
vertices. If $X$ is finite, the queue empties after finitely many steps, and every coarse vertex
outside $X$ adjacent to $X$ has been processed and declared bad. Apply Lemma~\ref{lem:contours} to
$X$: its boundary contains a cycle of some length $n$ from a family of at most $2n3^n$ candidates,
each with $k\ge n/7$ crossed edges with disjoint endpoints, each of which has a bad endpoint outside
$X$. Choosing one endpoint per edge in $2^k$ ways, Lemma~\ref{lem:badset} bounds the probability of a
given candidate by $(2f)^k\le(2f)^{n/7}$. By (S1),
$\Pp_\cI(|X|<\infty)\le\sum_n2n3^n(2f)^{n/7}<1$. Every good non-root vertex $v$ satisfies $H_v$, and
the boxes $B_v$ are disjoint, so with positive $\Pp_\cI$-probability the open cluster of $o$ in the
long-range model is infinite. Hence $\theta(q)\ge\Pp_q(\cI)\Pp_\cI(|X|=\infty)>0$ with $q<\beta_c$,
a contradiction.
\end{proof}

\section{The corollaries, and where the hypotheses enter}\label{sec:corollaries}

\begin{proof}[Proof of Corollary~\ref{cor:power}]
If $\alpha>d$ then $T(n)\le C'n^{-\alpha}$, so (Tail) holds, and (Pos) holds by the lower bound;
apply Theorem~\ref{thm:main}. If $\alpha<d$, put $s=d+\alpha<2d$ and $P_k=1-e^{-\beta J(k)}$ at a
parameter $\beta$ with $\theta(\beta)>0$. Berger's theorem~\cite[Theorem~2.6]{Berger2002} states: if
$\{P_k\}$ is percolating and $\liminf_{\norm{k}\to\infty}P_k\norm{k}_1^{s}>0$ for some $s<2d$, then
$\{(1-\epsilon)P_k\}$ is percolating for some $\epsilon>0$ (his standing assumption is
$0\le P_k=P_{-k}<1$, which holds here). Monotonicity and $1-e^{-t}\ge t/2$ on $[0,1]$ give
$P_k\ge(\beta c/2)\norm{k}^{-s}\ge(\beta c/2)\norm{k}_1^{-s}$ for large $\norm{k}$. Concavity of
$t\mapsto1-e^{-t}$ gives $1-e^{-(1-\epsilon)\beta J}\ge(1-\epsilon)(1-e^{-\beta J})$, so
$\theta((1-\epsilon)\beta)>0$. Thus $\theta(\beta)>0$ forces $\beta>\beta_c$.
\end{proof}

\begin{proof}[Proof of Corollary~\ref{cor:exp}]
Exponential decay and finite range both give $n^dT(n)\to0$.
\end{proof}

\begin{remark}[where the hypotheses enter]\label{rem:hyp}
The tail condition is used only in (S5), at one pair of scales $(w,r)$ with $r$ a large constant
multiple of $w$ fixed by (S6): in
Lemma~\ref{lem:faces}, for exits from $\Lam_\ell$ by edges longer than $w$, and in the event
$E_{\mathrm{long}}$ of Lemma~\ref{lem:extension}, for entrances that jump over the shell. The layer
argument needs $R\gg e^Lw$ and the relay needs $r\gg R+w$; both hold with $r\le C_0w$, and
(Tail) then makes $\beta(41r+1)^dT(w)\le\beta(42C_0)^dw^dT(w)$ small. Summability gives the uniform bound $a_z\le A$, which turns a large crossing
weight into many seed attempts. (Sym) gives the equal probabilities of the faces in
Lemma~\ref{lem:faces} and the orientation reductions. (Pos) gives $\beta_c<\infty$ and the
connectivity of open boxes. The dimension $d\ge2$ enters through the coordinate plane in
Lemma~\ref{lem:basic} and the coarse lattice $\Z^2$ with its contour count.
\end{remark}

\begin{remark}[the borderline]\label{rem:border}
For $J(x)\asymp\norm{x}^{-2d}$ the quantity $n^dT(n)$ stays bounded away from zero, and nothing in the
argument guarantees (S5) and (S6) together: long edges at scale $n$ leave the cube $[-n,n]^d$ with
probability of order one at every scale. In $d=1$ and $J(x)\asymp|x|^{-2}$ the transition is discontinuous
\cite{AizenmanNewman1986,DuminilCopinGarbanTassion2020}, and the argument above also lacks the coarse
plane. Apart from sufficiently spread-out kernels in $d\ge7$, where $d>6=3(\alpha\wedge2)$
and~\cite{HeydenreichHofstadSakai2008} applies, whether $\theta(\beta_c)=0$ for $\alpha=d\ge2$ remains open.
\end{remark}

{\small\raggedright

}

\end{document}